\documentclass[11pt]{article}
\usepackage{amsmath,amssymb,geometry,verbatim,enumerate,float,cite,setspace}
\newtheorem{theorem}{Theorem}
\newtheorem{lemma}[theorem]{Lemma}
\newtheorem{corollary}[theorem]{Corollary}
\usepackage{setspace}

\begin{document}

\onehalfspace

\title{\Large Graphs in which some and every maximum matching is uniquely restricted}

\author{Lucia Draque Penso$^1$, Dieter Rautenbach$^1$, U\'{e}verton dos Santos Souza$^2$}

\date{}

\maketitle

\begin{center}
{\small 
$^1$ Institute of Optimization and Operations Research, Ulm University, Ulm, Germany\\
\texttt{lucia.penso@uni-ulm.de, dieter.rautenbach@uni-ulm.de}\\[3mm]
$^2$
Instituto de Computa\c{c}\~{a}o, Universidade Federal Fluminense, Niter\' {o}i, Brazil\\
\normalsize \texttt{usouza@ic.uff.br}
}
\end{center}

\begin{abstract}
A matching $M$ in a graph $G$ is uniquely restricted 
if there is no matching $M'$ in $G$ 
that is distinct from $M$ 
but covers the same vertices as $M$.
Solving a problem posed by Golumbic, Hirst, and Lewenstein,
we characterize the graphs in which some maximum matching is uniquely restricted.
Solving a problem posed by Levit and Mandrescu,
we characterize the graphs in which every maximum matching is uniquely restricted.
Both our characterizations lead to efficient recognition algorithms for the corresponding graphs.
\end{abstract}

{\small 

\medskip


\medskip


}

\section{Introduction}\label{section1}

We consider finite and simple graphs as well as digraphs, and use standard terminology and notation.

A {\it matching} in a graph $G$ is a set of disjoint edges of $G$. 
A matching in $G$ of maximum cardinality is {\it maximum}.
A matching $M$ in $G$ is 
{\it perfect} if each vertex of $G$ is incident with an edge in $M$,
and 
{\it near-perfect} if each but exactly one vertex of $G$ is incident with an edge in $M$.
A graph $G$ is {\it factor-critical} if $G-u$ has a perfect matching for every vertex $u$ of $G$.
For a matching $M$ in $G$, let $V_G(M)$ denote the set of vertices of $G$ that are incident with an edge in $M$.
A path or cycle in $G$ is {\it $M$-alternating} if one of every two adjacent edges belongs to $M$.
For two sets $M$ and $N$, the {\it symmetric difference $M\Delta N$} 
is the set $(M\setminus N)\cup (N\setminus M)$.
Note that $\Delta$ is commutative and associative, 
that is, $M\Delta N=N\Delta M$ and $(M\Delta N)\Delta O=M\Delta (N\Delta O)$.
For a digraph $D$ and a vertex $u$ of $D$, 
let $V_D^+(u)$ be the set of vertices $v$ of $D$ 
such that $D$ contains a directed path from $u$ to $v$.
Similarly, let $V_D^-(u)$ be the set of vertices $w$ of $D$ 
such that $D$ contains a directed path from $w$ to $u$.
For a directed path or cycle $\vec{P}$, let $P$ denote the underlying undirected path or cycle.
For a positive integer $k$, let $[k]$ denote the set of positive integers at most $k$.
A set $I$ of vertices of a graph is {\it independent} 
if no two vertices in $I$ are adjacent.
An independent set of maximum cardinality is {\it maximum}.
Classical results of K\H{o}nig \cite{k} and Gallai \cite{g} imply that
$|I|+|M|=n$
for a bipartite graph $G$ of order $n$,
a maximum matching $M$ in $G$, and
a maximum independent set $I$ in $G$.

Golumbic, Hirst, and Lewenstein \cite{ghl} define a matching $M$ in a graph $G$ to be {\it uniquely restricted} 
if there is no matching $M'$ in $G$ with $M'\not=M$ and $V_G(M')=V_G(M)$,
that is, $M$ is the unique perfect matching in the subgraph $G[V_G(M)]$ of $G$ induced by $V_G(M)$.
In \cite{ghl} they show that it is NP-hard to determine a uniquely restricted matching of maximum size 
in a given bipartite graph that has a perfect matching.
Furthermore, they ask for which graphs 
the maximum size of a uniquely restricted matching
equals the size of a maximum matching, that is, 
for which graphs some maximum matching is uniquely restricted.
In \cite{lm} Levit and Mandrescu ask how to recognize the graphs
for which every maximum matching is uniquely restricted.
We answer both these questions completely
giving structural characterizations of both these classes of graphs
that lead to efficient recognition algorithms.

\section{Some maximum matching is uniquely restricted}\label{section2}

Let $I$ be an independent set in a bipartite graph $G$,
and let $\sigma:x_1,\ldots,x_k$ be a linear ordering of the elements of $I$.

For $j\in [k]$, let $I^\sigma_{\leq j}=\left\{ x_i:i\in[j]\right\}$. 

For $y\in N_G(I)$, let $p(y)=x_i$, where $i=\min\left\{ j\in [k]:y\in N_G\left(I^\sigma_{\leq j}\right)\right\}$,
that is, the index $i$ is such that $y\not\in N_G(x_1)\cup \ldots \cup N_G(x_{i-1})$ but $y\in N_G(x_i)$.
Let 
$$M^\sigma=\left\{ yp(y):y\in N_G(I)\right\}.$$
Note that in the graph $(V(G),M^\sigma)$, 
every vertex in $N_G(I)$ has degree exactly one.

If $E$ is a subset of the set $E(G)$ of edges of $G$, then $\sigma$ is {\it $E$-good} if $M^\sigma\subseteq E$.

The linear ordering $\sigma$ is an {\it accessibility ordering} for $I$ \cite{lm} if 
$$\left|N_G\left(I^\sigma_{\leq j}\right)\right|-\left|N_G\left(I^\sigma_{\leq j-1}\right)\right|\leq 1$$ 
for every $j\in [k]$.
Note that the definitions immediately imply that 
$\sigma$ is an accessibility ordering if and only if $M^\sigma$ is a matching in $G$.

A {\it partial accessibility ordering} for $I$ is an accessibility ordering $\sigma'$ for a subset $I'$ of $I$.

We summarize some results from \cite{ghl} that will be used.
\begin{theorem}[Golumbic, Hirst, and Lewenstein \cite{ghl}]\label{theorem-0}
A matching $M$ in a bipartite graph $G$ is uniquely restricted if and only if $G$ contains no $M$-alternating cycle.
\end{theorem}
The following result slightly extends Theorem 3.2 in \cite{lm}.

\begin{lemma}\label{theorem-1}
Let $G$ be a bipartite graph, and let $E$ be a set of edges of $G$.

The following statements are equivalent.
\begin{enumerate}[(i)]
\item There is a maximum independent set $I$ in $G$ 
that has an $E$-good accessibility ordering $\sigma$.
\item There is a maximum matching $M$ in $G$ such that $M$ is uniquely restricted and $M\subseteq E$.
\item Every maximum independent set $I$ in $G$ 
has an $E$-good accessibility ordering $\sigma$.
\end{enumerate}
\end{lemma}
{\it Proof:} (i) $\Rightarrow$ (ii). 
Let $I$ and $\sigma:x_1,\ldots,x_k$ be as in (i).
As noted above, $M^\sigma$ is a matching. 
Since $\sigma$ is $E$-good, we have $M^\sigma\subseteq E$.
By construction, $|M^\sigma|=|N_G(I)|$, and, 
since $I$ is a maximum independent set in $G$,
we have $|V(G)|=|I|+|N_G(I)|=|I|+|M^\sigma|$.
This implies that $M^\sigma$ is a maximum matching in $G$.
Let $\sigma':x'_1,\ldots,x'_\ell$ be the subordering of $\sigma$ formed by those $x_j$ 
where $j\in [k]$ is such that 
$\left|N_G\left(I^\sigma_{\leq j}\right)\right|-\left|N_G\left(I^\sigma_{\leq j-1}\right)\right|=1$,
that is, $\sigma'$ arises from $\sigma$ by removing the $x_j$ with 
$N_G(x_j)\subseteq N_G\left(I^\sigma_{\leq j-1}\right)$.
Let $N_G(I)=\{ y_1,\ldots,y_\ell\}$ be such that $p(y_i)=x_i'$ for $i\in [\ell]$,
that is, $M^\sigma=\{ x_i'y_i:i\in[\ell]\}$.
For a contradiction, we assume that $M^\sigma$ is not uniquely restricted.
By Theorem \ref{theorem-0},
there is an $M^\sigma$-alternating cycle $C$.
Since every edge of $C$ is incident with a vertex in $I$, and $I$ is independent, 
$C$ alternates between $I$ and $N_G(I)$, that is,
$C$ has the form $y_{r_1}x'_{r_1}y_{r_2}x'_{r_2}\ldots y_{r_t}x'_{r_t}y_{r_1}$.
Since $y_{r_i}\in N_G(x'_{r_{i-1}})\cap N_G(x'_{r_i})$ for $i\in [t]$, where we identify indices modulo $t$, 
the definition of $p(\cdot)$ implies the contradiction $r_1>r_2>r_3>\ldots >r_t>r_1$.
Hence, $M^\sigma$ is uniquely restricted,
and $G$ satisfies (ii).

\bigskip

\noindent (ii) $\Rightarrow$ (iii). 
Let $M=\{ x_1y_1,\ldots,x_\ell y_\ell\}$ be a maximum matching in $G$ 
such that $M$ is uniquely restricted and $M\subseteq E$. 
Let $I$ be a maximum independent set in $G$.
As noted in the introduction, we have $|I|+|M|=|V(G)|$.
Since $I$ contains at most one vertex from each edge in $M$,
this implies that $I$ contains all vertices in $V(G)\setminus V_G(M)$,
and exactly one vertex from each edge in $M$.
We may assume that $I=\{ x_1\ldots,x_\ell,x_{\ell+1}\ldots,x_k\}$,
where $V(G)\setminus V_G(M)=\{x_{\ell+1}\ldots,x_k\}$.
Note that the vertices $x_1,\ldots,x_\ell$ 
not necessarily belong to the same partite set of the bipartite graph $G$.
If there is some set $J\subseteq [\ell]$ 
such that 
$\left|N_G(x_j)\cap \left\{ y_i:i\in J\right\}\right|\geq 2$ for every $j\in J$,
then, since $I$ is independent, 
$G$ contains an $M$-alternating cycle, which is a contradiction.
Hence, for every set $J\subseteq [\ell]$, there is some $j\in J$ with  
$N_G(x_j)\cap \left\{ y_i:i\in J\right\}=\{ y_j\}$.
Therefore, we may assume that $x_1,\ldots,x_\ell$ 
are ordered in such a way that 
$i\geq j$ for every $i,j\in [\ell]$ with $x_iy_j\in E(G)$.
This implies that $\sigma:x_1\ldots x_k$ is an accessibility ordering for $I$
such that $M^\sigma=M\subseteq E$,
that is, $G$ satisfies (iii).

\bigskip

\noindent (iii) $\Rightarrow$ (i). This implication is trivial. $\Box$

\begin{lemma}\label{lemma-1}
Let $G$ be a bipartite graph, 
let $E$ be a set of edges of $G$,
and let $I$ be a maximum independent set in $G$.

$I$ has an $E$-good accessibility ordering $\sigma:x_1,\ldots,x_k$
if and only if 
for every $E$-good partial accessibility ordering $\sigma':x_1',\ldots,x'_{\ell-1}$ for $I$ with $0\leq \ell-1<|I|$,
there is an $E$-good partial accessibility ordering $\sigma'':x_1',\ldots,x'_{\ell-1},x'_\ell$ for $I$,
that is, every $E$-good partial accessibility ordering that does not contain all of $I$ can be extended.
\end{lemma}
{\it Proof:} 
Since the sufficiency is trivial, we only prove the necessity.
Let $\sigma$ and $\sigma'$ be as in the statement. 

If $\{ x_1,\ldots,x_{\ell-1}\}=\{ x'_1,\ldots,x'_{\ell-1}\}$, then
$N_G(x_\ell)\setminus N_G(\{ x'_1,\ldots,x'_{\ell-1}\})=N_G(x_\ell)\setminus N_G(\{ x_1,\ldots,x_{\ell-1}\})$.
Furthermore, 
if $N_G(x_\ell)\setminus N_G(\{ x'_1,\ldots,x'_{\ell-1}\})$ 
contains a vertex $y$,
then, since $\sigma$ is $E$-good, we have $x_\ell y\in E$.
Therefore, $\sigma'':x_1',\ldots,x'_{\ell-1},x_\ell$ is an $E$-good partial accessibility ordering for $I$.

If $\{ x_1,\ldots,x_{\ell-1}\}\not=\{ x'_1,\ldots,x'_{\ell-1}\}$, 
then $\{ x_1,\ldots,x_{\ell-1}\}\not\subseteq \{ x'_1,\ldots,x'_{\ell-1}\}$.
For $j=\min\{ i\in [\ell-1]:x_i\not\in \{ x'_1,\ldots,x'_{\ell-1}\}\}$,
we have $x_1,\ldots,x_{j-1}\in \{ x'_1,\ldots,x'_{\ell-1}\}$, and hence,
$N_G(x_j)\setminus N_G(\{ x'_1,\ldots,x'_{\ell-1}\})\subseteq 
N_G(x_j)\setminus N_G(\{ x_1,\ldots,x_{j-1}\})$.
Furthermore, 
if $N_G(x_j)\setminus N_G(\{ x'_1,\ldots,x'_{\ell-1}\})$ 
contains a vertex $y$,
then $y\in N_G(x_j)\setminus N_G(\{ x_1,\ldots,x_{j-1}\})$,
and hence, since $\sigma$ is $E$-good, we have $x_j y\in E$.
Therefore, $\sigma'':x_1',\ldots,x'_{\ell-1},x_j$ is an $E$-good partial accessibility ordering for $I$.
$\Box$

\begin{corollary}\label{corollary-1}
For a given bipartite graph $G$, and a given set $E$ of edges of $G$,
it is possible to check in polynomial time
whether $G$ has a maximum matching $M$
such that $M$ is uniquely restricted and $M\subseteq E$.
\end{corollary}
{\it Proof:} Since $G$ is bipartite, one can determine a maximum independent set $I$ in $G$ in polynomial time.
By Lemma \ref{theorem-1}, $G$ has the desired matching if and only if 
$I$ has an $E$-good accessibility ordering.
By Lemma \ref{lemma-1}, 
this can be checked by starting with the empty partial accessibility ordering for $I$, which is trivially $E$-good,
and iteratively extending $E$-good partial accessibility orderings for $I$ in a greedy way.
$\Box$

\bigskip

\noindent We now invoke the famous Gallai-Edmonds Structure Theorem \cite{lp},
which will be of central importance for this and the next section.

For a graph $G$,
\begin{itemize}
\item let $D(G)$ be the set of all vertices of $G$ that are not covered by some maximum matching in $G$,
\item let $A(G)$ be the set of vertices in $V(G)\setminus D(G)$ that have a neighbor in $D(G)$, and 
\item let $C(G)=V(G)\setminus (A(G)\cup D(G))$.
\end{itemize}
Let $G_B$ be the bipartite graph obtained from $G$ 
by deleting all vertices in $C(G)$ and all edges between vertices in $A(G)$, 
and by contracting each component $H$ of $G[D(G)]$ to a single vertex also denoted $H$.

Note that for a given graph $G$, the set $D(G)$,
and hence also $A(G)$ as well as $C(G)$,
can be determined in polynomial time \cite{lp}.

\begin{theorem}[Gallai-Edmonds Structure Theorem \cite{lp}]\label{theorem-4}
Let $G$ be a graph.

If $D(G)$, $A(G)$, $C(G)$, and $G_B$ are as above, then the following statements hold.
\begin{enumerate}[(i)]
\item Every component of $G[D(G)]$ is factor-critical.
\item Every component of $G[C(G)]$ has a perfect matching.
\item A matching in $G$ is maximum if and only if it is the union of 
\begin{enumerate}[(a)]
\item a near-perfect matching in each component of $G[D(G)]$,
\item a perfect matching in each component of $G[C(G)]$, and 
\item a matching with $|A(G)|$ edges
that matches the vertices in $A(G)$
with vertices in different components of $G[D(G)]$.
\end{enumerate}
\end{enumerate}
\end{theorem}
We proceed to the main result in this section.

\begin{theorem}\label{theorem-5}
Let $G$ be a graph.
Let $D(G)$, $A(G)$, $C(G)$, and $G_B$ be as above.
Let $E$ be the set of edges $aH$ of $G_B$, where $a\in A(G)$ and $H$ is a component of $G[D(G)]$,
such that the vertex $a$ has a unique neighbor, say $h$, in $V(H)$, 
and $H-h$ has a unique perfect matching.

Some maximum matching in $G$ is uniquely restricted if and only if
the following conditions hold.
\begin{enumerate}[(i)]
\item Every component of $G[C(G)]$ has a unique perfect matching.
\item $G_B$ has a maximum matching $M_B$ such that
\begin{enumerate}[(a)]
\item $M_B$ is uniquely restricted and 
\item $M_B\subseteq E$
\end{enumerate} 
\item Every component $H$ of $G[D(G)]$ has a vertex $h$ such that $H-h$ has a unique perfect matching.
\end{enumerate} 
\end{theorem} 
{\it Proof:} We first prove the necessity.
Therefore, let $M$ be a maximum matching in $G$ that is uniquely restricted.
Theorem \ref{theorem-4}(iii)(b) implies (i). 
Let $M_B$ be the matching in $G_B$ such that $M_B$ contains the edge $aH$,
where $a\in A$ and $H$ is a component of $G[D(G)]$,
if and only if $M$ contains an edge between the vertex $a$ and a vertex of $H$.
We will show that $M_B$ is as in (ii).
Theorem \ref{theorem-4}(iii)(c) implies that $M_B$ is a maximum matching of $G_B$.
If $M_B$ is not uniquely restricted, then
Theorem \ref{theorem-4}(i) and (iii) imply that $G$ has a maximum matching $M'$
with $V_G(M')=V_G(M)$ such that $M_B'\not=M_B$,
where $M_B'$ is defined analogously to $M_B$.
This implies $M'\not=M$, which is a contradiction.
Hence, (ii)(a) holds.
If some edge $aH$ in $M_B$ does not belong to $E$, 
then 
either $a$ has at least two distinct neighbors in $V(H)$
or $a$ has a unique neighbor $h$ in $V(H)$ but $H-h$ does not have a unique perfect matching.
In both cases,
Theorem \ref{theorem-4}(i) and (iii) imply that $G$ has a maximum matching $M'$
with $V_G(M')=V_G(M)$ that differs from $M$ within $H$,
which is a contradiction.
Hence, (ii)(b) holds.
If some component $H$ of $G[D(G)]$ has no vertex $h$ such that $H-h$ has a unique perfect matching, 
then Theorem \ref{theorem-4}(iii)(a) implies that $G$ has a maximum matching $M'$
with $V_G(M')=V_G(M)$ that differs from $M$ within $H$,
which is a contradiction.
Hence, (iii) holds.

Now we prove the sufficiency.
Let $M_1$ be the unique perfect matching in $G[C(G)]$.
Let $M_B$ be as in (ii).
Let $M_2$ be a matching in $G$ such that for every $a\in A$,
the matching $M_2$ contains an edge $ah$,
where $h\in V(H)$ and $H$ is a component of $G[D(G)]$, 
if and only if $M_B$ contains the edge $aH$.
By Theorem \ref{theorem-4}(iii)(c), $M_2$ covers all of $A(G)$.
By (ii)(b), $M_2$ is uniquely determined.
For every component $H$ of $G[D(G)]$ 
such that $M_2$ contains an edge $ah$ with $h\in V(H)$,
(ii)(b) implies that $H-h$ has a unique perfect matching $M_H$.
For every component $H$ of $G[D(G)]$ 
such that $M_2$ does not contain an edge $ah$ with $h\in V(H)$,
(iii) implies that $H$ has a vertex $h$ such that $H-h$ has a unique perfect matching $M_H$.
Let 
$$M_3=\bigcup_{H:H\,\,{\rm is\,\,a\,\,component\,\,of}\,\,G[D(G)]}M_H$$
and
$M=M_1\cup M_2\cup M_3$.
By Theorem \ref{theorem-4}(iii), $M$ is a maximum matching in $G$.
We will show that $M$ is uniquely restricted.
For a contradiction, we assume that $M'$ is a maximum matching in $G$ with $M'\not=M$ and $V_G(M')=V_G(M)$.
By (i) and Theorem \ref{theorem-4}(iii)(b), $M'$ contains $M_1$.
By (ii)(a) and (b), $M'$ contains $M_2$.
By (ii)(b) and (iii), $M'$ contains $M_3$.
Altogether, $M\subseteq M'$, which implies the contradiction $M=M'$. 
$\Box$

\begin{corollary}\label{corollary-2}
For a given graph $G$,
it is possible to check in polynomial time
whether some maximum matching in $G$ is uniquely restricted.
\end{corollary}
{\it Proof:} If some graph $H$ has a perfect matching $M$, 
then $M$ is uniquely restricted if and only if $H-e$ has no perfect matching for every $e\in M$.
Therefore, the conditions (i) and (iii) from Theorem \ref{theorem-5} can be checked in polynomial time.
By Corollary \ref{corollary-1}, condition (ii) from Theorem \ref{theorem-5} can be checked in polynomial time.
Now, Theorem \ref{theorem-5} implies the desired statement. $\Box$

\bigskip

\noindent Note that the constructive proofs of 
Lemma \ref{theorem-1},
Corollary \ref{corollary-1}, and
Theorem \ref{theorem-5}
also lead to an efficient algorithm that 
determines a maximum matching in a given graph $G$ 
that is uniquely restricted,
if such a matching exists.

\section{Every maximum matching is uniquely restricted}\label{section3}
It is convenient to split this section into two subsections, one about bipartite graphs, 
and one about not necessarily bipartite graphs.

\subsection{Bipartite graphs}\label{section3.1}

Throughout this subsection, let $G$ be a bipartite graph with partite sets $A$ and $B$.

For a matching $M$ in $G$,
let $D(M)$ be the digraph with vertex set $V(G)$ and arc set 
$$\{ (a,b):a\in A, b\in B,\mbox{ and }ab\in E(G)\setminus M\}\cup \{ (b,a):a\in A, b\in B,\mbox{ and }ab\in M\}.$$
Note that $M$-alternating paths and cycles in $G$ correspond to directed paths and cycles in $D$.

Let 
$$A_0(M)=\left\{ x\in A:d^-_{D(M)}(x)=0\right\}
\,\,\,\,\,\,\mbox{ and }\,\,\,\,\,\,
B_0(M)=\left\{ x\in B:d^+_{D(M)}(x)=0\right\}.$$
Note that 
\begin{eqnarray*}
A_0(M)&=&A\setminus V_G(M)\,\,\,\,\,\mbox{ and }\,\,\,\,\,\mbox{ $d^-_{D(M)}(a)=1$ for every $a\in A\setminus A_0(M)$},\\
B_0(M)&=&B\setminus V_G(M)\,\,\,\,\,\mbox{ and }\,\,\,\,\,\mbox{ $d^+_{D(M)}(b)=1$ for every $b\in B\setminus B_0(M)$.}
\end{eqnarray*}
Let 
$$V^+(M)=\bigcup_{a\in A_0(M)}V^+_{D(M)}(a)
\,\,\,\,\,\,\mbox{ and }\,\,\,\,\,\,
V^-(M)=\bigcup_{b\in B_0(M)}V^-_{D(M)}(b),$$
that is, 
$V^+(M)$ is the set of vertices of $G$ that are reachable from a vertex in $A_0(M)$ on an $M$-alternating path,
and
$V^-(M)$ is the set of vertices of $G$ that can reach a vertex in $B_0(M)$ on an $M$-alternating path.

\bigskip 

\noindent K\H{o}nig's classical method \cite{k} of finding a maximum matching in a bipartite graph 
relies on the following result (cf. Section 16.3 of \cite{s}).

\begin{theorem}[K\H{o}nig \cite{k}]\label{theorem1}
A matching $M$ in a bipartite graph $G$ is maximum 
if and only if $G$ contains no $M$-alternating path between a vertex in $A_0(M)$ and a vertex in $B_0(M)$,
that is, if and only if $V^+(M)\cap V^-(M)=\emptyset$.
\end{theorem}
In view of the correspondence between $M$-alternating cycles in $G$ and directed cycles in $D(M)$,
Golumbic, Hirst, and Lewenstein's \cite{ghl} 
characterization of a uniquely restricted matching in a bipartite graph can be rephrased as follows.
\begin{theorem}[Golumbic, Hirst, and Lewenstein \cite{ghl}]\label{theorem2}
A matching $M$ in a bipartite graph $G$ is uniquely restricted 
if and only if $D(M)$ is acyclic.
\end{theorem}
Our main result in this subsection is the following.

\begin{theorem}\label{theorem3}
Let $M$ be a maximum matching in a bipartite graph $G$.

Every maximum matching in $G$ is uniquely restricted if and only if 
$D(M)$ is acyclic,
and the two subgraphs $G[V^+(M)]$ and $G[V^-(M)]$ of $G$ induced by $V^+(M)$ and $V^-(M)$, respectively,
are forests.
\end{theorem}
The rest of this subsection is devoted to the proof of Theorem \ref{theorem3}.

\begin{lemma}\label{lemma1}
Let $M$ be a maximum matching in a bipartite graph $G$.

If $M'$ is a maximum matching in $G$, then $V^+(M')=V^+(M)$ and $V^-(M')=V^-(M)$.
\end{lemma}
{\it Proof:} Since the non-trivial components of $(V(G),M\Delta M')$
are $M$-$M'$-alternating cycles and
$M$-$M'$-alternating paths of even length,
it suffices, by an inductive argument, to show that 
$V^+(M')=V^+(M)$ and $V^-(M')=V^-(M)$
if 
either $M'=M\Delta E(C)$, 
where $C$ is an $M$-alternating cycle,
or $M'=M\Delta E(P)$, 
where $P$ is an $M$-alternating path between some vertex $a$ in $A_0(M)$ and some vertex $a'$ in $A\setminus A_0(M)$.
In the first case, $D(M')$ arises from $D(M)$ by inverting the orientation of the edges of $C$,
$A_0(M')=A_0(M)$, and $B_0(M')=B_0(M)$,
which easily implies $V^+(M')=V^+(M)$ and $V^-(M')=V^-(M)$.
Now, let $M'=M\Delta E(P)$, where $P$ is as above.
$D(M)$ contains a directed path $\vec{P}$ from $a$ to $a'$
such that $P$ is the underlying undirected path of $\vec{P}$.
Furthermore, $D(M')$ arises by inverting the orientation of the arcs of $\vec{P}$.
Since $M=M'\Delta E(P)$, $a'\in A_0(M')$, and $a\in A\setminus A_0(M')$, 
in order to complete the proof, it suffices, by symmetry, 
to show $V^+(M')\subseteq V^+(M)$ and $V^-(M')\subseteq V^-(M)$.

If $x\in V^+(M)\setminus V^+(M')$, 
then some directed path in $D(M)$ from a vertex in $A_0(M)$ to $x$ intersects $\vec{P}$,
which implies that $D(M')$ contains a directed path from $a'$ to $x$, that is, $x\in V_{D(M')}^+(a')\subseteq V^+(M')$,
which is a contradiction. Hence, $V^+(M')\subseteq V^+(M)$.
Similarly, if $x\in V^-(M)\setminus V^-(M')$, 
then some directed path in $D(M)$ from $x$ to a vertex $b$ in $B_0(M)$ intersects $\vec{P}$,
which implies that $D(M)$ contains a directed path from $a\in A_0(M)$ to $b\in B_0(M)$.
By Theorem \ref{theorem1}, $M$ is not maximum,
which is a contradiction.
$\Box$

\begin{lemma}\label{lemma2}
Let $M$ be a maximum matching in a bipartite graph $G$.

If every maximum matching in $G$ is uniquely restricted, 
then the two subgraphs $G[V^+(M)]$ and $G[V^-(M)]$ of $G$ induced by $V^+(M)$ and $V^-(M)$, respectively,
are forests.
\end{lemma}
{\it Proof:} For a contradiction, we may assume, by symmetry, that $G[V^+(M)]$ is not a forest.
For a cycle $C$ in $G[V^+(M)]$ and a maximum matching $M'$ in $G$, 
let $\vec{C}(M')$ be the subdigraph of $D(M')$ such that $C$ is the underlying undirected graph of $\vec{C}(M')$.
Since $M'$ is uniquely restricted, Theorem \ref{theorem2} implies that $\vec{C}(M')$ is not a directed cycle in $D(M')$.
Therefore, the set
$$S(C,M')=\left\{ x\in V(C):d^-_{\vec{C}(M')}(x)=0\right\}$$
is not empty.
Note that $\left|\left\{ x\in V(C):d^+_{\vec{C}(M')}(x)=0\right\}\right|=|S(C,M')|$,
that is, $\vec{C}(M')$ contains equally many sink vertices as source vertices.

We assume that $C$ and $M'$ are chosen such that $|S(C,M')|$ is minimum.

Let $x\in S(C,M')$. 
Since $d^+_{\vec{C}(M')}(x)=2$, we have $x\in A$.
Since $x\in V(C)\subseteq V^+(M)$, 
Lemma \ref{lemma1} implies $x\in V^+(M')$.
Hence, there is a directed path $\vec{P}$ in $D(M')$ from some vertex $a$ in $A_0(M')$ to $x$.
First, we assume that $\vec{P}$ and $\vec{C}(M')$ only share the vertex $x$.
Let $\vec{Q}$ be a directed path in $\vec{C}(M')$ from $x$ to some vertex $y\in V(C)$ with $d^+_{\vec{C}(M')}(y)=0$.
Since $d^-_{\vec{C}(M')}(y)=2$, we have $y\in B$.
Since $y\in V_{D(M')}^+(a)\subseteq V^+(M')$, Theorem \ref{theorem1} implies $y\in B\setminus B_0(M')$.
This implies that there is some vertex $a'$ such that $a'y\in M'$.
If $\vec{R}$ is the concatenation of $\vec{P}$, $\vec{Q}$, and the arc $(y,a')$, 
and $M''=M'\Delta E(R)$,
then $|S(C,M'')|$ is strictly smaller than $|S(C,M')|$, 
which is a contradiction.
Hence, $\vec{P}$ and $\vec{C}(M')$ share a vertex different from $x$.
This implies that $\vec{P}$ contains a directed subpath $\vec{P}'$ 
from a vertex $y$ in $V(C)\setminus \{ x\}$ to $x$ 
such that $\vec{P}'$ is internally disjoint from $\vec{C}(M')$.
If $z$ is such that $(z,y)$ is an arc of $\vec{C}(M')$,
and $Q$ is the path in $C$ between $x$ and $y$ that contains $z$,
then $C'=P'\cup Q$ is a cycle in $G[V^+(M)]$ such that 
$|S(C',M')|$ is strictly smaller than $|S(C,M')|$,
which is a contradiction.
Hence, we may assume that $d^-_{\vec{C}(M')}(y)=0$.
Now, if $R$ is one of the two paths in $C$ between $x$ and $y$,
then $C''=P'\cup R$ is a cycle in $G[V^+(M)]$ such that 
$|S(C'',M')|$ is strictly smaller than $|S(C,M')|$,
which is a contradiction.
$\Box$

\bigskip

\noindent If $M$ is a maximum matching in $G$, and
$a\in A_0(M)$ and $a'b'\in M$ are such that $b'$ is a neighbor of $a$, 
then $M'=(M\setminus \{ a'b'\})\cup \{ ab'\}$ is a maximum matching in $G$,
and we say that $M'$ arises from $M$ by an {\it edge exchange}.
Similarly, 
if $b\in B_0(M)$ and $a'b'\in M$ are such that $a'$ is a neighbor of $b$, 
then $M''=(M\setminus \{ a'b'\})\cup \{ a'b\}$ is a maximum matching in $G$,
and also in this case, we say that $M''$ arises from $M$ by an edge exchange.

\begin{lemma}\label{lemma3}
Let $M$ be a maximum matching in a bipartite graph $G$.

If $D(M)$ is acyclic, then every maximum matching in $G$ 
arises from $M$ by a sequence of edge exchanges.
\end{lemma}
{\it Proof:} If $M'$ is any maximum matching in $G$, then, since $D$ is acyclic, 
the non-trivial components of $(V(G),M\Delta M')$ 
are $M$-$M'$-alternating paths $P_1,\ldots,P_k$,
each starting with an edge in $M$ and ending with an edge in $M'$.
Clearly, $M'=M\Delta E(P_1)\Delta\cdots \Delta E(P_k)$.
Since the maximum matching $M\Delta E(P_1)$ arises from $M$ by a sequence of edge exchanges,
the statement follows easily by an inductive argument. $\Box$

\begin{lemma}\label{lemma4}
Let $M$ be a maximum matching in a bipartite graph $G$.
Let $D(M)$ be acyclic, 
and let the two subgraphs $G[V^+(M)]$ and $G[V^-(M)]$ of $G$ induced by $V^+(M)$ and $V^-(M)$, respectively,
be forests.

If $M'$ arises from $M$ by an edge exchange,
then $D(M')$ is acyclic.
\end{lemma}
{\it Proof:} By symmetry, we may assume that 
$a\in A_0(M)$ and $a'b'\in M$ are such that $b'$ is a neighbor of $a$, 
and that $M'=(M\setminus \{ a'b'\})\cup \{ ab'\}$.
Note that $A_0(M')=(A_0(M)\setminus \{ a\})\cup \{ a'\}$,
$d^-_{D(M)}(a)=0$, and $d^-_{D(M')}(a')=0$.
If $\vec{C}$ is a directed cycle in $D(M')$, 
then, since $D(M)$ is acyclic, $\vec{C}$ contains the arc $(b',a)$ of $D(M')$.
This implies that $G[V^+_{D(M)}(a)]$, and hence, also $G[V^+(M)]$ contains the cycle $C$,
which is a contradiction.
Hence, $D(M')$ is acyclic.
$\Box$

\bigskip

\noindent We are now in a position to prove Theorem \ref{theorem3}.

\bigskip

\noindent {\it Proof of Theorem \ref{theorem3}:}
The necessity follows from Theorem \ref{theorem2} and Lemma \ref{lemma2}.
For the sufficiency, let $M'$ be any maximum matching of $G$.
By Lemma \ref{lemma3}, $M'$ arises from $M$ by a sequence of edge exchanges.
By Lemma \ref{lemma1} and Lemma \ref{lemma4}, 
it follows by induction on the number of these edge exchanges
that $D(M')$ is acyclic.
Therefore, by Theorem \ref{theorem2}, $M'$ is uniquely restricted. $\Box$

\subsection{Not necessarily bipartite graphs}\label{section3.2}

In order to extend Theorem \ref{theorem3} to graphs that are not necessarily bipartite,
we again rely on the Gallai-Edmonds Structure Theorem.

\begin{theorem}\label{theorem5}
Let $G$ be a graph. 
Let $D(G)$, $A(G)$, $C(G)$, and $G_B$ be as above.

Every maximum matching in $G$ is uniquely restricted if and only if 
the following conditions hold.
\begin{enumerate}[(i)]
\item Every component of $G[C(G)]$ has a unique perfect matching.
\item For every component $H$ of $G[D(G)]$, 
every near-perfect matching in $H$ is uniquely restricted.
\item Every maximum matching of $G_B$ is uniquely restricted.
\item If an edge $aH$ of $G_B$, where $a\in A(G)$ and $H$ is a component of $G[D(G)]$,
is contained in some maximum matching of $G_B$, 
then the vertex $a$ has a unique neighbor in $V(H)$.
\end{enumerate} 
\end{theorem}
{\it Proof:} In view of Theorem \ref{theorem-4}(iii), the proof of the necessity is straightforward;
in fact, it can be done using very similar arguments as the proof of the necessity in Theorem \ref{theorem-5}.
Therefore, we proceed to show the sufficiency.
Let $M$ be a maximum matching in $G$.
By Theorem \ref{theorem-4}(iii)(b), (i) implies that $M\cap E(G[C(G)])$ is uniquely determined.
By Theorem \ref{theorem-4}(iii)(a) and (c), (iii) and (iv) imply that 
$M\cap \{ uv\in E(G):u\in A(G)\mbox{ and }v\in D(G)\}$ is uniquely determined,
which also implies that for every component $H$ of $G[D(G)]$, 
the unique vertex of $H$ that is not covered by an edge in 
$M\cap E(G[D(G)])$ is uniquely determined.
Now, by Theorem \ref{theorem-4}(iii)(a), (ii) implies that 
$M\cap E(G[D(G)])$ is uniquely determined,
which completes the proof. $\Box$

\bigskip

\noindent Note that the factor-critical graphs in which every near-perfect matching is uniquely restricted
(cf. Theorem \ref{theorem5}(ii))
are exactly the factor-critical graphs $G$ with the minimum possible number $|V(G)|$ of distinct near-perfect matchings.
In \cite{dr} it is shown that these are exactly the connected graphs whose blocks are odd cycles. 

\begin{corollary}\label{corollary1}
The graphs $G$ with the property that every maximum matching in $G$ is uniquely restricted 
can be recognized in polynomial time.
\end{corollary}
{\it Proof:} Theorem \ref{theorem3} obviously implies the statement if $G$ is bipartite.
As noted above the sets $D(G)$, $A(G)$, and $C(G)$ can be determined in polynomial time for a given graph $G$.
If $G$ has a perfect matching $M$, then $M$ is uniquely restricted
if and only if $G-e$ has no perfect matching for every $e\in M$.
If $G$ has a near-perfect matching $M$ that does not cover the vertex $u$ of $G$, 
then $M$ is uniquely restricted if and only if 
$M$ is a uniquely restricted perfect matching of $G-u$.
Since it is easy to check in polynomial time 
whether some edge of a bipartite graph belongs to some maximum matching,
and also 
whether some vertex of a bipartite graph is not covered by some maximum matching,
the four conditions in Theorem \ref{theorem5} can be checked in polynomial time,
which completes the proof. $\Box$


\begin{thebibliography}{}
\bibitem{dr} T. Do\v{s}li\'{c} and D. Rautenbach, Factor-Critical Graphs with the Minimum Number of Near-Perfect Matchings, manuscript 2015.
\bibitem{g} T. Gallai, \"{U}ber extreme Punkt- und Kantenmengen, Ann. Univ. Sci. Budapest. E\"{o}tv\"{o}s Sect. Math. 2 (1959) 133-138.
\bibitem{ghl} M.C. Golumbic, T. Hirst, and M. Lewenstein, Uniquely restricted matchings, Algorithmica 31 (2001) 139-154.
\bibitem{k} D. K\H{o}nig, Graphok \'{e}s matrixok, Matematikai \'{e}s Fizikai Lapok 38 (1931) 116-119.
\bibitem{lm} V.E. Levit and E. Mandrescu, Local maximum stable sets in bipartite graphs with uniquely restricted maximum matchings, Discrete Appl. Math. 132 (2003) 163-174
\bibitem{lp} L. Lov\'{a}sz and M. Plummer, Matching Theory, North-Holland, 1986.
\bibitem{s} A. Schrijver, Combinatorial Optimization - Polyhedra and Efficiency, Springer 2004.
\end{thebibliography}
\end{document}